\documentclass{article}

\usepackage{amssymb,exscale,amsmath}

\author{Katharina Habermann\footnote{SUB G\"ottingen,
    Platz der G\"ottinger Sieben 1, 37073 G\"ottingen, Germany,
    {\it e-mail: habermann@sub.uni-goettingen.de}},
  Lutz Habermann\footnote{Department of Mathematics,
    University of Hannover, Welfengarten 1, 30167 Hannover, Germany,
    {\it e-mail: habermann@math.uni-hannover.de}},
  Paul Rosenthal\footnote{Department of Mathematics and Computer Science,
    University of Greifswald, Jahnstr. 15a, 17487 Greifswald, Germany,
    {\it e-mail: paul.rosenthal@uni-greifswald.de}}}

\title{Symplectic Yang-Mills theory, Ricci tensor, and connections
  \footnote{This work was partially supported by the Deutsche
  Forschungsgemeinschaft (DFG).}}

\date{\empty}

\numberwithin{equation}{section}

\newcommand{\vp}{\varphi}
\newcommand{\rd}{\mathrm{d}}
\newcommand{\R}{\mathbb{R}}
\newcommand{\cC}{\mathcal{C}}
\newcommand{\cE}{\mathcal{E}}
\newcommand{\mfe}{\mathfrak{e}}
\newcommand{\mfs}{\mathfrak{s}}
\newcommand{\sfb}{\mathsf{b}}
\newcommand{\sfe}{\mathsf{e}}
\newcommand{\sfs}{\mathsf{s}}
\newcommand{\Tr}{\operatorname{Tr}}
\newcommand{\Sec}{{\mathrm{\Gamma}}}
\newcommand{\End}{{\mathrm{End}}}
\newcommand{\ric}{{\mathrm{ric}}}
\newcommand{\sric}{{\mathrm{sric}}}
\newcommand{\sRic}{{\mathrm{sRic}}}

\newtheorem{lm}{Lemma}[section]
\newtheorem{sz}[lm]{Proposition}
\newtheorem{tm}[lm]{Theorem}
\newtheorem{fo}[lm]{Corollary}

\newtheorem{bm0}[lm]{Remark}
\newtheorem{bs0}[lm]{Example}

\newenvironment{beweis}{{\em Proof.}}{\hspace*{\fill} $\square$}
\newenvironment{bm}{\begin{bm0}\rm}{\hspace*{\fill} $\square$ \end{bm0}}

\begin{document}

\maketitle

\begin{abstract}
  A Yang--Mills theory in a purely symplectic framework is developed.
  The corresponding Euler--Lagrange equations are derived and first
  integrals are given. We relate the results to the work of
  Bourgeois and Cahen on preferred symplectic connections.
\end{abstract}

\section{Introduction}

The study of connections in vector resp. principal fiber bundles over
manifolds has been a broad field of research in many contexts over decades. In
particular, one is interested in properties of connections respecting various
types of geometries.

In Riemannian geometry, we have the classical result that
each Riemannian manifold $(M,g)$ admits a distinguished linear connection,
called the Levi--Civita connection. This connection is characterized
by the conditions that it is metric, i.e. $g$ is parallel, and that its
torsion vanishes. More general, for any $2$-form $\zeta$ on $M$ with
values in the tangent bundle $TM$
of $M$, there exists a unique metric connection on $(M,g)$ having $\zeta$
as its torsion.

Turning one's interest to symplectic geometry, the situation changes
drastically. Namely, if $(M,\omega)$ is an almost symplectic manifold,
then, for a given $2$-form $\zeta$ on $M$ with values in $TM$,
the space of symplectic connections whose torsion is $\zeta$
is either empty or an infinite dimensional affine space
(see Lemma~\ref{lm:torsion-sympl}). Here, by a symplectic connection,
we mean the analog of a metric connection on a Riemannian manifold.
Accordingly, a symplectic connection on $(M,\omega)$ is a linear
connection $\nabla$ on $M$ such that $\nabla\omega=0$.

By the above, on an almost symplectic manifold, there is no canonical
symplectic connection singled out by any torsion condition, i.e.
there is no notion analogous to the Levi--Civita connection.
This fact gives rise to the question whether it is possible to select a single
class of symplectic connections in other ways.
In recent years, this problem has been approached by
studying symplectic connections that, in addition to torsion-freeness,
satisfy suitable curvature conditions.
For an overview, see \cite{bielgutt} and the references
therein. In \cite{bourcah}, the authors derived such a condition
by a variational principle using a Lagrangian density with a
quadratic polynomial in the curvature.
The connections satisfying the deduced Euler--Lagrange equations are now
referred to as preferred symplectic connections.
Unfortunately, aside from the surface case, there is not much known
about general properties of those connections. Progress has been made
towards understanding so-called symplectic connections of Ricci type,
which are described by a curvature condition that implies the
Euler--Lagrange equations of \cite{bourcah}. But all these studies
do not consider possible effects caused by the torsion. It is
generally assumed that all connections taken into consideration
are torsion-free.

It is the aim of the present work to develop an approach that also
includes symplectic connections with non-vanishing torsion.
Treating symplectic connections without any torsion obstruction has several
advantages and is motivated by the following aspects.

Also in Riemannian geometry, in spite of the existence of the
Levi--Civita connection, one is more and more interested in connections
with torsion. This is due to current issues in string theory. There
one studies different types of metric connections having good but
non-trivial torsion. Here, in most cases, ``good'' means that the
torsion considered as a covariant $3$-tensor field is
totally skew-symmetric (see e.g. \cite{agricfried}).
Another point is that dropping the restriction to torsion-free connections
allows to take into account connections that are
Hermitian with respect to a compatible almost complex structure. 
In general, those connections have non-trivial torsion.
On the other hand, for fixed compatible almost complex structure,
there exists a distinguished Hermitian connection (see \cite{gauduchon}).
A further reason comes from the theory of symplectic Dirac operators.
The more one studies these operators, the more one becomes convinced
that, at least in this context, certain symplectic connections with torsion
are more suitable than torsion-free connections (see \cite{unser-buch}).
Moreover, only admitting symplectic connections with torsion
makes it possible to extend the considerations to
almost symplectic manifolds (cp. Proposition~\ref{sz:ex-torsion-free}).

The approach given here has the advantage that it also works for
connections in vector bundles. Generalizing the ansatz in \cite{bourcah},
it gives a purely symplectic Yang--Mills theory.
At this point, we should mention the paper \cite{urakawa}.
At first glance, this paper is closely related to that presented here,
and really there are some relations to our work. However, Urakawa
studies the usual Yang--Mills functional, using the symplectic
framework. More precisely, he considers the variation of the integral
\begin{displaymath}
  \int_M g\left(R^\nabla,R^\nabla\right)\omega^n\;,
\end{displaymath}
where $g$ is a compatible Riemannian metric on the symplectic
manifold $(M,\omega)$.

Dealing with Yang--Mills theory on an almost symplectic manifold, the
so-called symplectic Ricci tensor comes naturally into play.
This tensor field generalizes the usual Ricci tensor for
a torsion-free symplectic connection as well as the
$\ast$-Ricci tensor of an almost Hermitian manifold. Furthermore,
it has the advantage that the corresponding endomorphism,
called the symplectic Ricci operator, can be defined also for
connections in a vector bundle.

The paper is organized as follows. In the second and third section,
we recall some results on the torsion of metric and symplectic connections
and on Ricci tensors, respectively. In Section~4, we introduce
two functionals. The first one is the symplectic analog of the
Yang--Mills functional, whereas the second one is defined by the
symplectic Ricci operator. We compute their
Euler--Lagrange equations and describe how they are related.
In Section~5, we discuss first integrals of the symplectic
Yang--Mills equation, which, in dimension 4, are symplectic
analogs of the self-duality and anti-self-duality equation.
In the last section, we relate our results to that in
\cite{bourcah}. In particular, we give a proof for the
Euler--Lagrange equations of the functionals considered there.

\section{Connections and torsion}

In this section, we summarize some well known facts
on the torsion of metric and symplectic connections.

Let $M$ be a smooth manifold. If $E$ is a vector bundle over $M$, we denote
the space of smooth sections of $E$ by $\Sec(E)$, the space of
smooth $k$-forms on $M$ with values in $E$ by $\Omega^k(M,E)$, and the
endomorphism bundle of $E$ by $\End(E)$.

Let $\cC(M)$ be the space of linear connections on $M$, i.e. of connections
in the tangent bundle $TM$ of $M$. Then $\cC(M)$ is an affine space
over the vector space $\Omega^1(M,\End(TM))$. The torsion
$T^\nabla\in\Omega^2(M,TM)$ of a connection $\nabla\in\cC(M)$ is
defined by
\begin{displaymath}
  T^\nabla(X,Y)=\nabla_X Y-\nabla_Y X-[X,Y]
\end{displaymath}
for $X,Y\in\Sec(TM)$. Thus the map
$\nabla\in\cC(M)\mapsto T^\nabla\in\Omega^2(M,TM)$
is affine and the associated linear map
$\Phi:\Omega^1(M,\End(TM))\to\Omega^2(M,TM)$ is given by
\begin{displaymath}
  \Phi(\theta)(X,Y)=\theta(X)Y-\theta(Y)X\;.
\end{displaymath}
A connection $\nabla\in\cC(M)$ is called torsion-free if $T^\nabla=0$.

Let $g$ be a Riemannian metric on $M$. A metric connection on $M$
is a connection $\nabla\in\cC(M)$ such that $\nabla g=0$, i.e.
\begin{displaymath}
  X(g(Y,Z))=g(\nabla_X Y,Z)+g(Y,\nabla_X Z)
\end{displaymath}
for all $X,Y,Z\in\Sec(TM)$. The space $\cC(M,g)$ of metric connections
on $M$ is an affine subspace of $\cC(M)$ and its corresponding vector
space $\cE^1(M,g)$ is formed by all $\theta\in\Omega^1(M,\End(TM))$ that
satisfy
\begin{equation}
  \label{eq:diff-metr-zsh}
  g(\theta(X)Y,Z)=-g(\theta(X)Z,Y)\;.
\end{equation}

\begin{lm}
  The restriction of $\Phi$ to $\cE^1(M,g)$ is an isomorphism onto
  the space $\Omega^2(M,TM)$.
\end{lm}
\begin{beweis}
  Let $\theta\in\cE^1(M,g)$ and suppose that $\Phi(\theta)=0$, i.e.
  \begin{displaymath}
    \theta(X)Y=\theta(Y)X\;.
  \end{displaymath}
  By means of Equation~(\ref{eq:diff-metr-zsh}), we conclude
  \begin{align*}
    g(\theta(X)Y,Z)
    &=-g(\theta(X)Z,Y)
    =-g(\theta(Z)X,Y)
    =g(\theta(Z)Y,X)\\
    &=g(\theta(Y)Z,X)
    =-g(\theta(Y)X,Z)
    =-g(\theta(X)Y,Z)\;,
  \end{align*}
  which implies $\theta=0$.

  Let $\zeta\in\Omega^2(M,TM)$ and let $\theta\in\Omega^1(M,\End(TM))$
  be determined by
  \begin{displaymath}
    2\,g(\theta(X)Y,Z)
    =g(\zeta(X,Y),Z)+g(\zeta(Z,X),Y)+g(\zeta(Z,Y),X)\;.
  \end{displaymath}
  One readily verifies that $\theta\in\cE^1(M,g)$ and that
  \begin{displaymath}
    g(\theta(X)Y,Z)-g(\theta(Y)X,Z)
    =g(\zeta(X,Y),Z)\;,
  \end{displaymath}
  which means $\Phi(\theta)=\zeta$.
\end{beweis}

\medskip

An immediate consequence is

\begin{fo}
  The map $\nabla\in\cC(M,g)\mapsto T^\nabla\in\Omega^2(M,TM)$
  is $1:1$. In particular, there is a unique connection
  $\nabla\in\cC(M,g)$, called the Levi--Civita connection, such
  that $T^\nabla=0$.
  \hspace*{\fill} $\square$
\end{fo}

From now on, we suppose that $M$ is endowed with an almost symplectic
structure, i.e. a non-degenerate $2$-form $\omega$. The form $\omega$
is called a symplectic structure if in addition it is closed. A connection
$\nabla\in\cC(M)$ is said to be symplectic if $\nabla\omega=0$, i.e. if
\begin{displaymath}
  X(\omega(Y,Z))=\omega(\nabla_X Y,Z)+\omega(Y,\nabla_X Z)
\end{displaymath}
for all $X,Y,Z\in\Sec(TM)$. We point out that we do not require that
a symplectic connection is torsion-free. The space $\cC(M,\omega)$ of
symplectic connections on $M$ is again an affine subspace of $\cC(M)$
and its vector space $\cE^1(M,\omega)$ consists of all
$\theta\in\Omega^1(M,\End(TM))$ that
satisfy
\begin{displaymath}
  \omega(\theta(X)Y,Z)=\omega(\theta(X)Z,Y)\;.
\end{displaymath}
Although a symplectic connection is the analog of a metric connection,
the properties of the torsion maps are completely different
(cf. \cite {tondeur}).

\begin{lm}
  \label{lm:torsion-sympl}
  \begin{itemize}
    \item[{\rm (i)}] A form $\theta\in\cE^1(M,\omega)$ satisfies
      $\Phi(\theta)=0$ if and only if the expression $\omega(\theta(X)Y,Z)$
      is totally symmetric in $X,Y,Z\in\Sec(TM)$.
    \item[{\rm (ii)}] The image $\Phi\left(\cE^1(M,\omega)\right)$
      of $\cE^1(M,\omega)$ under $\Phi$ is the space of all
      $\zeta\in\Omega^2(M,TM)$ such that
      \begin{equation}
        \label{eq:image-phi-om}
        \omega(\zeta(X,Y),Z)+\omega(\zeta(Y,Z),X)+\omega(\zeta(Z,X),Y)=0\;.
      \end{equation}
    \end{itemize}
\end{lm}
\begin{beweis}
  Assertion~(i) is obvious. Moreover, it is easy to check that
  Equation~(\ref{eq:image-phi-om}) holds true for any
  $\zeta\in\Phi\left(\cE^1(M,\omega)\right)$. Finally, suppose that
  $\zeta\in\Omega^2(M,TM)$ satisfies Equation~(\ref{eq:image-phi-om})
  and let $\theta\in\Omega^1(M,\End(TM))$ be given by
  \begin{displaymath}
    \omega(\theta(X)Y,Z)
    =\frac{1}{3}(\omega(\zeta(X,Y),Z)+\omega(\zeta(X,Z),Y))\;.
  \end{displaymath}
  Then $\theta\in\cE^1(M,\omega)$ and $\Phi(\theta)=\zeta$.
\end{beweis}

\begin{fo}
  The map $\nabla\in\cC(M,\omega)\mapsto T^\nabla\in\Omega^2(M,TM)$
  is neither injective nor onto. Moreover, the pre-image
  $\left\{\nabla\in\cC(M,\omega):T^\nabla=\zeta\right\}$ of a form
  $\zeta\in\Omega^2(M,TM)$ is either empty or infinite dimensional.
  \hspace*{\fill} $\square$
\end{fo}

A further difference to the metric case is

\begin{sz}
  \label{sz:ex-torsion-free}
  There exists a torsion-free symplectic connection $\nabla$ on $M$
  if and only if $\omega$ is a symplectic structure.
\end{sz}
\begin{beweis}
  This follows from the following two facts. For any
  $\nabla\in\cC(M,\omega)$, we have
  \begin{displaymath}
    \rd\omega(X,Y,Z)=\omega\left(T^\nabla(X,Y),Z\right)
    +\omega\left(T^\nabla(Y,Z),X\right)
    +\omega\left(T^\nabla(Z,X),Y\right)\;.
  \end{displaymath}
  If $\nabla^0\in\cC(M)$ is any torsion-free connection and
  $\theta\in\Omega^1(M,\End(TM))$ is defined by
  \begin{displaymath}
    \omega(\theta(X)Y,Z)
    =\frac{1}{3}\left(\left(\nabla^0_X\omega\right)(Y,Z)
        +\left(\nabla^0_Y\omega\right)(X,Z)\right)\;,
  \end{displaymath}
  then the connection $\nabla=\nabla^0+\theta$ is also torsion-free
  and
  \begin{displaymath}
    X(\omega(Y,Z))-\omega(\nabla_X Y,Z)-\omega(Y,\nabla_X Z)
    =\frac{1}{3}\,\rd\omega(X,Y,Z)\;.
  \end{displaymath}
\end{beweis}

\medskip

According to the above, in symplectic geometry, there is no analog of the
Levi-Civita connection. Furthermore, in the case that $\omega$ is a
symplectic structure, any connection $\nabla\in\cC(M,\omega)$ such that
$\omega\left(T^\nabla(X,Y),Z\right)$ is symmetric or skew-symmetric
in $Y,Z$ has to be torsion-free.

\section{Ricci tensors}

In the following, let $M$ have dimension $2n$ and let
$\sfs=(\sfe_1,\ldots,\sfe_{2n})$ be a
symplectic frame on some open subset $U\subset M$, i.e. a frame
of vector fields on $U$ such that
\begin{displaymath}
  \omega(\sfe_i,\sfe_j)=\omega(\sfe_{n+i},\sfe_{n+j})=0
  \quad\mbox{and}\quad
  \omega(\sfe_i,\sfe_{n+j})=\delta_{ij}
\end{displaymath}
for $i,j=1,\ldots,n$. Furthermore, let $J^\sfs$ be the almost complex
structure on $U$ defined by
\begin{displaymath}
  J^\sfs\sfe_i=\sfe_{n+i}
\end{displaymath}
for $i=1,\ldots,n$.

Let $\nabla$ be any connection on $M$. Its curvature
is the form $R^\nabla\in\Omega^2(M,\End(TM))$ given by
\begin{equation}
  \label{eq:curvature}
  R^\nabla(X,Y)Z=\nabla_X\nabla_Y Z-\nabla_Y\nabla_X Z
  -\nabla_{[X,Y]} Z\;.
\end{equation}
The Ricci tensor of $\nabla$ is the tensor field
$\ric^\nabla\in\Sec(T^\ast M\otimes T^\ast M)$ defined by
\begin{align*}
  \ric^\nabla(X,Y)
  &=\Tr\left(Z\mapsto R^\nabla(Z,X)Y\right)\\
  &=\sum_{i=1}^{2n}\omega\left(R^\nabla(\sfe_i,X)Y,J^\sfs\sfe_i\right)\;.
\end{align*}
Thus the Ricci tensor $\ric^\nabla$ is obtained by contracting the
curvature $R^\nabla$ with respect to the symplectic form $\omega$.
Following the ideas of Vaisman \cite{vaisman}, we consider another
contraction of $R^\nabla$. We define $\sRic^\nabla\in\Sec(\End(TM))$ by
\begin{displaymath}
  \sRic^\nabla(X)=\sum_{i=1}^n R^\nabla(\sfe_i,J^\sfs\sfe_i)X
\end{displaymath}
and call it the symplectic Ricci operator. Moreover, we define
the symplectic Ricci tensor as the tensor field
$\sric^\nabla\in\Sec(T^\ast M\otimes T^\ast M)$ given by
\begin{displaymath}
  \sric^\nabla(X,Y)
  =\omega\left(\sRic^\nabla(X),Y\right)\;.
\end{displaymath}

\begin{bm}
  If $\omega$ is the K\"ahler form of an almost Hermitian
  structure $(g,J)$ and $\nabla$ is the Levi-Civita connection
  of $g$, then $\sric^\nabla$ is the so-called $\ast$-Ricci tensor.
  See e.g. \cite{yano}.
\end{bm}

The Ricci tensors $\ric^\nabla$ and $\sric^\nabla$ are related by

\begin{sz}
  \label{sz:comp-ric-sric}
  If $\nabla\in\cC(M,\omega)$, then
  \begin{displaymath}
    \sric^\nabla(X,Y)-\ric^\nabla(X,Y)
    =\omega\left(K^\nabla(X),Y\right)\;,
  \end{displaymath}
  where
  \begin{align*}
    K^\nabla(X)
    &=\sum_{i=1}^n\left(
      T^\nabla\left(T^\nabla(\sfe_i,J^\sfs\sfe_i),X\right)
      +\left(\nabla_X T^\nabla\right)(\sfe_i,J^\sfs\sfe_i)\right)\\
    &\hspace{5mm}+\sum_{i=1}^{2n}\left(
      T^\nabla\left(T^\nabla(X,\sfe_i),J^\sfs\sfe_i\right)
      +\left(\nabla_{\sfe_i}T^\nabla\right)(J^\sfs\sfe_i,X)\right)\;.
  \end{align*}
\end{sz}
\begin{beweis}
  Let $\nabla\in\cC(M,\omega)$. Then
  \begin{equation}
    \label{eq:prop-symp-curv}
    \omega\left(R^\nabla(X,Y)Z_1,Z_2\right)
    =\omega\left(R^\nabla(X,Y)Z_2,Z_1\right)\;.
  \end{equation}
  With this, one gets
  \begin{align*}
    \lefteqn{\sric^\nabla(X,Y)-\ric^\nabla(X,Y)}\\
    &=\sum_{i=1}^n\omega\left(R^\nabla(\sfe_i,J^\sfs\sfe_i)X
      +R^\nabla(J^\sfs\sfe_i,X)\sfe_i
      +R^\nabla(X,\sfe_i)J^\sfs\sfe_i,Y\right)\;.
  \end{align*}
  Now, applying the first Bianchi identity
  (cf. \cite{kobnom1}, Chapter~III, Theorem~5.3),
  the assertion follows.
\end{beweis}

\medskip

Consequently, for any torsion-free symplectic connection $\nabla$, the
Ricci tensors $\ric^\nabla$ and $\sric^\nabla$ coincide.
But in general, for a generic symplectic connection, this is not true.

\begin{sz}
  If $\nabla\in\cC(M,\omega)$, then
  \begin{displaymath}
    \sric^\nabla(X,Y)=\sric^\nabla(Y,X)
  \end{displaymath}
  for any $X,Y\in\Sec(TM)$.
\end{sz}
\begin{beweis}
  This easily follows from Equation~(\ref{eq:prop-symp-curv}).
\end{beweis}

\medskip

By Proposition~\ref{sz:comp-ric-sric}, in general, $\ric^\nabla$ is not
symmetric. This, among others, indicates that the symplectic Ricci
tensor $\sric^\nabla$ is more adapted to symplectic geometry than
the usual Ricci tensor $\ric^\nabla$. Furthermore, the symplectic
Ricci operator $\sRic^\nabla$ can be straightforwardly generalized to
connections $\nabla$ in a vector bundle on $M$.
If $E$ is a vector bundle on $M$ and
$\nabla:\Sec(E)\to\Omega^1(M,E)$ is a connection in $E$, we define
the symplectic Ricci operator $\sRic^\nabla\in\Sec(\End(E))$
of $\nabla$ by
\begin{displaymath}
  \sRic^\nabla(\xi)=\sum_{i=1}^n R^\nabla(\sfe_i,J^\sfs\sfe_i)\xi
\end{displaymath}
for a section $\xi\in\Sec(E)$, where the curvature
$R^\nabla\in\Omega^2(M,\End(E))$ is given
analogously to Equation~(\ref{eq:curvature}) by
\begin{displaymath}
  R^\nabla(X,Y)\xi=\nabla_X\nabla_Y\xi-\nabla_Y\nabla_X\xi
  -\nabla_{[X,Y]}\xi\;.
\end{displaymath}

\section{Symplectic Yang--Mills functionals}

In this section, we generalize the variational principle for
symplectic connections suggested by Bourgeois and Cahen \cite{bourcah}
to connections in vector bundles.

We now suppose that the manifold $M$ is closed.
Let $E$ be a real vector bundle over $M$ of rank $2m$ with an almost
symplectic structure, i.e. a non-degenerate $2$-form
$\sfb\in\Sec\left(\Lambda^2 E\right)$. Let
$\mfs=(\mfe_1,\ldots,\mfe_{2m})$ be a local symplectic frame in $E$
and let $J^\mfs$ be the local almost complex structure in $E$ given by
\begin{displaymath}
  J^\mfs\mfe_i=\mfe_{m+i}
  \quad\mbox{for}\quad i=1,\ldots,m\;.
\end{displaymath}
We define pairings
\begin{displaymath}
  (K,L)\in\Sec(\End(E))\times\Sec(\End(E))
  \mapsto\sfb(K,L)\in C^\infty(M)
\end{displaymath}
and
\begin{displaymath}
  (\alpha,\beta)\in\Omega^k(M,\End(E))\times\Omega^k(M,\End(E))
  \mapsto\sfb(\alpha,\beta)\in C^\infty(M)
\end{displaymath}
by
\begin{displaymath}
  \sfb(K,L)=\sum_{i=1}^{2m}\sfb(K\mfe_i,LJ^\mfs\mfe_i)
\end{displaymath}
and
\begin{displaymath}
  \sfb(\alpha,\beta)=\sum_{1\le i_1<\cdots <i_k\le 2n}
  \sfb(\alpha(\sfe_{i_1},\ldots,\sfe_{i_k}),
  \beta(J^\sfs\sfe_{i_1},\ldots,J^\sfs\sfe_{i_k}))\;.
\end{displaymath}
It is easy to see that $\sfb(K,L)$ and $\sfb(\alpha,\beta)$ do not
depend on the choice of the symplectic frames $\sfs$ and $\mfs$
and that
\begin{displaymath}
  \sfb(K,L)=\sfb(L,K)
  \quad\mbox{and}\quad
  \sfb(\alpha,\beta)=(-1)^k\sfb(\beta,\alpha)\;.
\end{displaymath}

Let $\cC(E)$ denote the space of connections in $E$ and let
$\cC(E,\sfb)$ be the subspace of symplectic connections in $E$,
i.e. of connections $\nabla\in\cC(E)$ such that
\begin{displaymath}
  X(\sfb(\xi_1,\xi_2))=\sfb(\nabla_X\xi_1,\xi_2)+\sfb(\xi_2,\nabla_X\xi_1)
\end{displaymath}
for all $X\in\Sec(TM)$ and $\xi_1,\xi_2\in\Sec(E)$. Then $\cC(E)$
is an affine space over the vector space $\Omega^1(M,\End(E))$
and $\cC(E,\sfb)$ is an affine subspace over the vector
space $\cE^1(E,\sfb)$ of all $\theta\in\Omega^1(M,\End(E))$ that satisfy
\begin{displaymath}
  \sfb(\theta(X)\xi_1,\xi_2)=\sfb(\theta(X)\xi_2,\xi_1)\;.
\end{displaymath}

The symplectic analog of the Yang--Mills functional is now the
functional
\begin{displaymath}
  I_1:\cC(E,\sfb)\to\R\;,\qquad
  I_1(\nabla)=\frac{1}{2}\int_M
  \sfb\left(R^\nabla,R^\nabla\right)\omega^{(n)}\;.
\end{displaymath}
Besides this, we consider the functional
\begin{displaymath}
  I_2:\cC(E,\sfb)\to\R\;,\qquad
  I_2(\nabla)=\frac{1}{2}\int_M
  \sfb\left(\sRic^\nabla,\sRic^\nabla\right)\omega^{(n)}\;.
\end{displaymath}
Here we have used the abbreviation
\begin{displaymath}
  \omega^{(k)}=\frac{1}{k!}\omega^k\;.
\end{displaymath}
In particular, $\omega^{(n)}$ is the symplectic volume form.

To compute the Euler--Lagrange equations of the functionals $I_1$
and $I_2$, we need some preparations. Let $\Omega^k(M)$ denote
the space of smooth $k$-forms on $M$ and let
\begin{displaymath}
  (\alpha,\alpha')\in\Omega^k(M,\End(E))\times\Omega^l(M,\End(E))
  \mapsto\sfb(\alpha\wedge\alpha')\in\Omega^{k+l}(M)
\end{displaymath}
be the bilinear map determined by
\begin{displaymath}
  \sfb((K\otimes\vp)\wedge(L\otimes\psi))
  =\sfb(K,L)\vp\wedge\psi
\end{displaymath}
for $K,L\in\Sec(\End(E))$, $\vp\in\Omega^k(M)$ and
$\psi\in\Omega^l(M)$. The symplectic analog of the Riemannian Hodge
operator is defined as follows (cf. \cite{brylinski}).
The symplectic Hodge operator is the unique isomorphism
$\ast:\Omega^k(M,\End(E))\to\Omega^l(M,\End(E))$ that
satisfies
\begin{displaymath}
  \sfb(\alpha,\beta)\,\omega^{(n)}
  =\sfb(\alpha\wedge\ast\beta)
\end{displaymath}
for any $\alpha,\beta\in\Omega^k(M,\End(E))$. In the next lemma,
we note some properties of this operator (cf. \cite{frihaber}).

\begin{lm}
  \label{lm:prop-hodge}
  \begin{itemize}
  \item[{\rm (i)}] For any $\alpha\in\Omega^k(M,\End(E))$ and any $k$,
    it is $\ast(\ast\alpha)=\alpha$.
  \item[{\rm (iii)}] For any $L\in\Sec(\End(E))$, it is
    $\ast L=L\otimes\omega^{(n)}$.
  \item[{\rm (iii)}] If $\alpha\in\Omega^1(M,\End(E))$, then
    $\ast\alpha=\alpha\wedge\omega^{(n-1)}$.
  \item[{\rm (iv)}] If $\alpha\in\Omega^2(M,\End(E))$, then
    \begin{displaymath}
      \ast\alpha=\alpha(\omega)\otimes\omega^{(n-1)}
      -\alpha\wedge\omega^{(n-2)}\;.
    \end{displaymath}
    \hspace*{\fill} $\square$
  \end{itemize}
\end{lm}

Here $\alpha(\omega)\in\Sec(\End(E))$ for
$\alpha\in\Omega^2(M,\End(E))$ means the section given by
\begin{displaymath}
  \alpha(\omega)
  =\sum_{i=1}^n\alpha(\sfe_i,J^\sfs\sfe_i)\;.
\end{displaymath}

If $\nabla$ is a connection in $E$, let
$\rd^\nabla:\Omega^k(M,\End(E))\to\Omega^{k+1}(M,\End(E))$
be the associated exterior differential and set
\begin{displaymath}
  \delta^\nabla=(-1)^{k+1}{\ast\rd^\nabla\ast}:
  \Omega^{k+1}(M,\End(E))\to\Omega^k(M,\End(E))\;.
\end{displaymath}

Then $\delta^\nabla$ is the formal adjoint of $\rd^\nabla$ in the
following sense.

\begin{sz}
  \label{sz:adjoint}
  Let $\alpha\in\Omega^k(M,\End(E))$ and $\beta\in\Omega^{k+1}(M,\End(E))$.
  Then, for any $\nabla\in\cC(E,\sfb)$,
  \begin{displaymath}
    \int_M\sfb\left(\rd^\nabla\alpha,\beta\right)\omega^{(n)}
    =\int_M\sfb\left(\alpha,\delta^\nabla\beta\right)\omega^{(n)}\;.
  \end{displaymath}
\end{sz}
\begin{beweis}
  One can proceed as in the Riemannian case. Since the connection $\nabla$
  is symplectic, we have
  \begin{displaymath}
    \rd(\sfb(\alpha\wedge{\ast\beta}))
    =\sfb\left(\rd^\nabla\alpha\wedge{\ast\beta}\right)
      +(-1)^k\sfb\left(\alpha\wedge\rd^\nabla{\ast\beta}\right)\;.
  \end{displaymath}
  Using in addition Stokes' Theorem and Lemma~\ref{lm:prop-hodge}(i),
  we conclude
  \begin{align*}
    \int_M\sfb\left(\rd^\nabla\alpha,\beta\right)\omega^{(n)}
    &=\int_M\left(\rd^\nabla\alpha\wedge{\ast\beta}\right)\omega^{(n)}\\
    &=\int_M\rd(\sfb(\alpha\wedge{\ast\beta}))\omega^{(n)}
    -(-1)^k\int_M\sfb\left(\alpha\wedge\rd^\nabla{\ast\beta}\right)
    \omega^{(n)}\\
    &=(-1)^{k+1}\int_M\sfb\left(\alpha
      \wedge{\ast{\ast\rd^\nabla}}{\ast\beta}\right)\omega^{(n)}\\
    &=\int_M\sfb\left(\alpha,\delta^\nabla\beta\right)\omega^{(n)}\;.
  \end{align*}
\end{beweis}

\medskip

For the action of $\delta^\nabla$ on $1$-forms, we have
(cp. also \cite{brylinski}, \S 1.2 and Theorem~2.2.1)

\begin{lm}
  \label{lm:delta-on-1-forms}
  Suppose that $\omega$ is a symplectic structure and let $\nabla$ be
  any connection in $E$. Then
  \begin{displaymath}
    \delta^\nabla\alpha=-\rd^\nabla\alpha(\omega)
  \end{displaymath}
  for $\alpha\in\Omega^1(M,\End(E))$.
\end{lm}
\begin{beweis}
  By means of Lemma~\ref{lm:prop-hodge}, we derive
  \begin{align*}
    \delta^\nabla
    &=-{\ast\rd^\nabla}{\ast\alpha}\\
    &=-{\ast\rd^\nabla}\left(\alpha\wedge\omega^{(n-1)}\right)\\
    &={-\ast}\left(\rd^\nabla\alpha\wedge\omega^{(n-1)}\right)\\
    &={-\ast}\left(\rd^\nabla\alpha(\omega)\otimes\omega^{(n)}\right)\\
    &=-\rd^\nabla\alpha(\omega)\;.
  \end{align*}
\end{beweis}

\medskip

Now we can prove

\begin{tm}
  \label{tm:euler-lagrange}
  \begin{itemize}
  \item[{\rm (i)}] A connection $\nabla\in\cC(E,\sfb)$ is a critical
    point of the functional $I_1$ if and only if
    $\rd^\nabla{\ast R^\nabla}=0$.
  \item[{\rm (ii)}] In case $\omega$ is a symplectic structure,
    a connection $\nabla\in\cC(E,\sfb)$ is a critical point of the
    functional $I_2$ if and only if $\nabla\sRic^\nabla=0$. Moreover,
    in this case, the Euler--Lagrange equations of $I_1$ and $I_2$
    are equivalent.
  \end{itemize}
\end{tm}
\begin{beweis}
  Let $\theta\in\cE^1(E,\sfb)$ and let $\nabla^t$ be a smooth curve
  in $\cC(E,\sfb)$ with $\nabla^0=\nabla$ and
  \begin{displaymath}
    \left.\frac{\rd}{\rd t}\nabla^t\right|_{t=0}=\theta\;.
  \end{displaymath}
  Then, as is well known,
  \begin{displaymath}
    \left.\frac{\rd}{\rd t}R^{\nabla^t}\right|_{t=0}
    =\rd^\nabla\theta\;.
  \end{displaymath}
  Thus, by means of Proposition~\ref{sz:adjoint}, we get
  \begin{displaymath}
    \left.\frac{\rd}{\rd t}I_1\left(\nabla^t\right)\right|_{t=0}
    =\int_M\sfb\left(\rd^\nabla\theta,R^\nabla\right)\omega^{(n)}
    =\int_M\sfb\left(\theta,\delta^\nabla R^\nabla\right)\omega^{(n)}\;.
  \end{displaymath}
  Since $\delta^\nabla R^\nabla\in\cE^1(E,\sfb)$ because of the
  symplecticity of $\nabla$ and since also the restriction of
  $\sfb$ to $\cE^1(E,\sfb)$ is non-degenerate, Assertion~(i) follows.

  Now we suppose that $\omega$ is a symplectic structure. Then we
  can apply Lemma~\ref{lm:delta-on-1-forms} to obtain
  \begin{displaymath}
    \left.\frac{\rd}{\rd t}\mathrm{sRic}^{\nabla^t}\right|_{t=0}
    =\left.\frac{\rd}{\rd t}R^{\nabla^t}(\omega)\right|_{t=0}
    =\rd^\nabla\theta(\omega)=-\delta^\nabla\theta\;.
  \end{displaymath}
  Hence
  \begin{displaymath}
    \left.\frac{\rd}{\rd t}I_2\left(\nabla^t\right)\right|_{t=0}
    =-\int_M\sfb\left(\delta^\nabla\theta,\sRic^\nabla\right)\omega^{(n)}
    =-\int_M\sfb\left(\theta,\nabla\sRic^\nabla\right)\omega^{(n)}\;.
  \end{displaymath}
  Since also $\nabla\sRic^\nabla\in\cE^1(E,\sfb)$, this gives the first
  part of Assertion~(ii). The second part is a consequence of
  Proposition~\ref{sz:diff-I1-I2} below. Alternatively, it follows from
  \begin{align*}
    \rd^\nabla{\ast R^\nabla}
    &=\rd^\nabla\left(R^\nabla(\omega)\otimes\omega^{(n-1)}
      -R^\nabla\wedge\omega^{(n-2)}\right)\\
    &=\rd^\nabla\left(\sRic^\nabla\otimes\omega^{(n-1)}\right)
    -\rd^\nabla R^\nabla\wedge\omega^{(n-2)}
    -R^\nabla\wedge\rd\omega^{(n-2)}\\
    &=\nabla\sRic^\nabla\wedge\omega^{(n-1)}\\
    &=\ast\nabla\sRic^\nabla\;,
  \end{align*}
  where we have used again Lemma~\ref{lm:prop-hodge} and the
  Bianchi identity $\rd^\nabla R^\nabla=0$.
\end{beweis}

\begin{sz}
  \label{sz:diff-I1-I2}
  If $\omega$ is a symplectic structure, then the functionals
  $I_1$ and $I_2$ differ by a constant.
\end{sz}
\begin{beweis}
  We can proceed as in \cite{bourcah}, \S 2. First we observe that,
  by Lemma~\ref{lm:prop-hodge}(iv),
  \begin{align*}
    \sfb(\alpha\wedge\beta)\wedge\omega^{(n-2)}
    &=\sfb\left(\alpha\wedge\left(\beta\wedge\omega^{(n-2)}\right)\right)\\
    &=\sfb\left(\alpha\wedge\left(\beta(\omega)
        \otimes\omega^{(n-1)}\right)\right)
    -\sfb(\alpha\wedge\ast\beta)\\
    &=(\sfb(\alpha(\omega),\beta(\omega))-\sfb(\alpha,\beta))\,\omega^{(n)}
  \end{align*}
  for all $\alpha,\beta\in\Omega^2(M,\End(E))$. In particular,
  \begin{displaymath}
    \left(\sfb\left(\sRic^\nabla,\sRic^\nabla\right)
      -\sfb\left(R^\nabla,R^\nabla\right)\right)\omega^{(n)}
    =\sfb\left(R^\nabla\wedge R^\nabla\right)\wedge\omega^{(n-2)}\;.
  \end{displaymath}
  Hence, by $\rd^\nabla R^\nabla=0$, $\rd\omega=0$ and
  Stokes' Theorem, for $\nabla^t$ as in the proof of
  Theorem~\ref{tm:euler-lagrange}, we have
  \begin{align*}
    \left.\frac{\rd}{\rd t}(I_2-I_1)\left(\nabla^t\right)\right|_{t=0}
    &=\int_M\sfb\left(\rd^\nabla\theta\wedge R^\nabla\right)
    \wedge\omega^{(n-2)}\\
    &=\int_M\rd\left(\sfb\left(\theta\wedge R^\nabla\right)
      \wedge\omega^{(n-2)}\right)\\
    &=0\;.
  \end{align*}
\end{beweis}

\begin{bm}
  \begin{itemize}
  \item[{\rm (i)}] One can proceed analogously if $E$ is endowed with
    a Riemannian structure or if $E$ is a complex vector bundle with
    a Hermitian structure.
  \item[{\rm (ii)}] In the considered situation,
    $\sRic^\nabla=\lambda\operatorname{id}_E$ for some $\lambda\in\R$
    implies $\sRic^\nabla=0$. However, if $E$ is a Hermitian
    vector bundle, it may be interesting to study the analog of
    the Hermitian Yang--Mills equation.
  \end{itemize}
\end{bm}

\section{First integrals}

In the case that $M$ is $4$-dimensional, one can ask for self-dual
and anti-self-dual solutions of the symplectic Yang--Mills equation
\begin{equation}
  \label{eq:exmpl-ym}
  \rd^\nabla{\ast R^\nabla}=0\;,
\end{equation}
i.e. for connections $\nabla\in\cC(E,\sfb)$ such that
\begin{equation}
  \label{eq:self-dual}
  \ast R^\nabla=R^\nabla
\end{equation}
and
\begin{equation}
  \label{eq:anti-self-dual}
  \ast R^\nabla=-R^\nabla\;,
\end{equation}
respectively.

\begin{sz}
  \label{sz:self-dual}
  If $\dim M=4$, then
  \begin{itemize}
  \item[{\rm (i)}] $\ast R^\nabla=R^\nabla$ if and only
    \begin{displaymath}
      R^\nabla=\frac{1}{2}\,\sRic^\nabla\otimes\omega\;.
    \end{displaymath}
  \item[{\rm (ii)}] $\ast R^\nabla=-R^\nabla$ if and only if
    $\sRic^\nabla=0$.
  \end{itemize}
\end{sz}
\begin{beweis}
  This follows from
  \begin{displaymath}
    \ast R^\nabla=\sRic^\nabla\otimes\omega-R^\nabla\;.
  \end{displaymath}
\end{beweis}

\medskip

The self-duality equation~(\ref{eq:self-dual}) as well as the
anti-self-duality equation~(\ref{eq:anti-self-dual}) possess
generalizations to arbitrary dimensions. For the second equation,
this is obvious. Namely, by Proposition~\ref{sz:self-dual},
Equation~(\ref{eq:anti-self-dual}) can be generalized by
$\sRic^\nabla=0$. Concerning the self-duality equation, observe that
\begin{equation}
  \label{eq:gen-self-dual}
  R^\nabla=\frac{1}{n}\,\sRic^\nabla\otimes\omega
\end{equation}
is equivalent to the existence of an endomorphism $L\in\Sec(\End(E))$
such that $R^\nabla=L\otimes\omega$. Therefore,
Equation~(\ref{eq:gen-self-dual}) generalizes Equation~(\ref{eq:self-dual}).

\begin{lm}
  \label{lm:gen-sd}
  \begin{itemize}
  \item[{\rm (i)}] If $n=1$, then
    $R^\nabla=\sRic^\nabla\otimes\omega$ for any connection
    $\nabla\in\cC(E)$.
  \item[{\rm (ii)}] If $n\ge 2$, then Equation~(\ref{eq:gen-self-dual})
    is equivalent to
    \begin{equation}
      \label{eq:gen-sd-equ-1}
      \ast R^\nabla=\frac{1}{n-1}\,R^\nabla\wedge\omega^{(n-2)}\;.
    \end{equation}
  \item[{\rm (iii)}] It is $\sRic^\nabla=0$ if and only if
    \begin{equation}
      \label{eq:gen-sd-equ-2}
      \ast R^\nabla=-R^\nabla\wedge\omega^{(n-2)}\;.
    \end{equation}
  \end{itemize}
\end{lm}
\begin{beweis}
  The first assertion is trivial. To see the second assertion,
  suppose that $n\ge 2$ and assume first that
  Equation~(\ref{eq:gen-sd-equ-1}) holds true. Since
  \begin{equation}
    \label{eq:proof-gen-sd}
    \ast R^\nabla=\sRic^\nabla\otimes\omega^{(n-1)}
    -R^\nabla\wedge\omega^{(n-2)}
  \end{equation}
  by Lemma~\ref{lm:prop-hodge}(iv), this implies
  \begin{displaymath}
    \sRic^\nabla\otimes\omega^{(n-1)}=n\,{\ast R^\nabla}\;,
  \end{displaymath}
  which is equivalent to Equation~(\ref{eq:gen-self-dual}).
  The converse can be derived straightforwardly.
  Assertion~(iii) is an easy consequence of
  Equation~(\ref{eq:proof-gen-sd}).
\end{beweis}

\begin{fo}
  In case $\omega$ is a symplectic structure, any connection
  $\nabla\in\cC(E)$ such that $\sRic^\nabla=0$ is a solution
  of Equation~(\ref{eq:exmpl-ym}). If, in addition, $n\ge 2$,
  the same holds true for any connection $\nabla\in\cC(E)$
  that solves Equation~(\ref{eq:gen-self-dual}).
\end{fo}
\begin{beweis}
  This follows from the Bianchi identity $\rd^\nabla R^\nabla=0$
  and Lemma~\ref{lm:gen-sd}.
\end{beweis}

\section{Preferred symplectic connections}

The aim of this section is to relate the considerations of the previous two
sections to the original variational principle suggested by Bourgeois
and Cahen. For this, we suppose that $\omega$ is a symplectic structure
on $M$. Let $\cC_0(M,\omega)$ denote the space of torsion-free symplectic
connections. According to Lemma~\ref{lm:torsion-sympl} and
Proposition~\ref{sz:ex-torsion-free}, $\cC_0(M,\omega)$ is
an affine space over the vector space $\cE_0^1(M,\omega)$ of all
$\theta\in\Omega^1(M,\End(TM))$ such that $\omega(\theta(X)Y,Z)$
is totally symmetric in $X,Y,Z\in\Sec(TM)$.
Let $I_0:\cC_0(M,\omega)\to\R$ be the restriction
of the functional $I_2$ for the case $E=TM$ to the subspace $\cC_0(M,\omega)$.
The critical points of the functional $I_0$ are called
preferred symplectic connections (cf. \cite{bielgutt,cagura}).

\begin{tm}\cite{bourcah}
  \label{tm:bourcah}
  A connection $\nabla\in\cC_0(M,\omega)$ is a critical point
  of the functional $I_0$ if and only if
  \begin{displaymath}
    \left(\nabla_X\sric^\nabla\right)(Y,Z)
    +\left(\nabla_Y\sric^\nabla\right)(Z,X)
    +\left(\nabla_Z\sric^\nabla\right)(X,Y)=0
  \end{displaymath}
  for all $X,Y,Z\in\Sec(TM)$.
\end{tm}

For completeness and since, as it seems to us, there is no proof of this
result in the literature, we will give a proof of it here.

Let $\sigma:\Omega^1(M,\End(TM))\to\Omega^1(M,\End(TM))$ be the
Bianchi projector. That means that $\sigma(\theta)$ for
$\theta\in\Omega^1(M,\End(TM))$ is given by
\begin{displaymath}
  \omega(\sigma(\theta)(X)Y,Z)
  =\frac{1}{3}\,(\omega(\theta(X)Y,Z)+\omega(\theta(Y)Z,X)
  +\omega(\theta(Z)X,Y))\;.
\end{displaymath}
One easily checks

\begin{lm}
  \label{lm:bianchi}
  For any $\theta_1,\theta_2\in\Omega^1(M,\End(TM))$,
  \begin{displaymath}
    \omega(\sigma(\theta_1),\theta_2)
    =\omega(\theta_1,\sigma(\theta_2))\;.
  \end{displaymath}
  \hspace*{\fill} $\square$
\end{lm}

Furthermore, we have

\begin{lm}
  \label{lm:tang-tor-free}
  It is $\cE^1_0(M,\omega)=\sigma\left(\cE^1(M,\omega)\right)$.
  \hspace*{\fill} $\square$
\end{lm}
\begin{beweis}
  A direct calculation shows that
  $\sigma\left(\cE^1(M,\omega)\right)\subset\cE^1(M,\omega)$.
  Therefore the assertion is an easy consequence of $\sigma^2=\sigma$ and
  \begin{displaymath}
    \cE^1_0(M,\omega)=\left\{\theta\in\cE^1(M,\omega):
      \sigma(\theta)=\theta\right\}\;.
  \end{displaymath}
\end{beweis}

\medskip

{\em Proof of Theorem~\ref{tm:bourcah}.} According to the proof of
Theorem~\ref{tm:euler-lagrange}, a connection $\nabla\in\cC_0(M,\omega)$
is a critical point of $I_0$ if and only if
\begin{displaymath}
  \int_M\omega\left(\theta,\nabla\sRic^\nabla\right)\omega^{(n)}=0
\end{displaymath}
for all $\theta\in\cE^1_0(M,\omega)$. By Lemma~\ref{lm:tang-tor-free},
the last condition is equivalent to
\begin{displaymath}
  \int_M\omega\left(\sigma(\theta),\nabla\sRic^\nabla\right)\omega^{(n)}=0
\end{displaymath}
for all $\theta\in\cE^1(M,\omega)$. By Lemma~\ref{lm:bianchi}, this is
the same as
\begin{displaymath}
  \int_M\omega\left(\theta,\sigma\left(\nabla\sRic^\nabla\right)\right)
  \omega^{(n)}=0
\end{displaymath}
for all $\theta\in\cE^1(M,\omega)$. Since
$\nabla\sRic^\nabla\in\cE^1(M,\omega)$, we obtain that
$\nabla\in\cC_0(M,\omega)$ is a critical point of $I_0$ if and only if
\begin{displaymath}
  \sigma\left(\nabla\sRic^\nabla\right)=0\;,
\end{displaymath}
which, because of
\begin{displaymath}
  \omega\left(\left(\nabla_X\sRic^\nabla\right)Y,Z\right)
  =\left(\nabla_X\sric^\nabla\right)(Y,Z)\;,
\end{displaymath}
is the desired relation.
\hspace*{\fill} $\square$

\begin{bm}
  The condition for a connection $\nabla\in\cC_0(M,\omega)$ to be
  preferred can also be expressed as
  \begin{displaymath}
    \sigma\left(\delta^\nabla R^\nabla\right)=0\;,
  \end{displaymath}
  since $\nabla\sRic^\nabla=\delta^\nabla R^\nabla$ as shown in the
  proof of Theorem~\ref{tm:euler-lagrange}.
\end{bm}

We conclude with the following observation.

\begin{sz}
  Let $n\ge 2$. Then any torsion-free connection $\nabla\in\cC(M)$
  that solves Equation~(\ref{eq:gen-self-dual}) has to be flat.
\end{sz}
\begin{beweis}
  Let $\nabla\in\cC(M)$ be a torsion-free solution of
  Equation~(\ref{eq:gen-self-dual}). Then, by the first Bianchi identity,
  we have
  \begin{align*}
    0&=
    \sum_{i=1}^n\left(R^\nabla(\sfe_i,J^\sfs\sfe_i)X
      +R^\nabla(J^\sfs\sfe_i,X)\sfe_i+R^\nabla(X,\sfe_i)J^\sfs\sfe_i\right)\\
    &=\frac{1}{n}\,\sum_{i=1}^n\left(
      \omega(\sfe_i,J^\sfs\sfe_i)\sRic^\nabla(X)
      +\omega(J^\sfs\sfe_i,X)\sRic^\nabla(\sfe_i)
      +\omega(X,\sfe_i)\sRic^\nabla(J^\sfs\sfe_i)\right)\\
    &=\frac{1}{n}\left(n\,\sRic^\nabla(X)
      +\sRic^\nabla\left(
        \sum_{i=1}^n(\omega(J^\sfs\sfe_i,X)\sfe_i-\omega(\sfe_i,X)J^\sfs\sfe_i)
      \right)\right)\\
    &=\frac{n-1}{n}\,\sRic^\nabla(X)\;,
  \end{align*}
  which yields $\sRic^\nabla=0$. Thus, by Lemma~\ref{lm:gen-sd},
  the connection $\nabla$ satisfies Equation~(\ref{eq:gen-sd-equ-1})
  as well as Equation~(\ref{eq:gen-sd-equ-2}), and this implies
  $R^\nabla=0$.
\end{beweis}

\end{document}